\documentclass[conference]{IEEEtran}

\usepackage{amsmath,amssymb,amsthm}
\usepackage{graphicx}
\usepackage{bm,color}
\usepackage{algorithm,algpseudocode}
\usepackage{url}
\usepackage{cite}
\usepackage{textcomp}

\allowdisplaybreaks[2]

\newtheorem{proposition}{Proposition}
\theoremstyle{definition}
\newtheorem{remark}{Remark}

\newcommand{\cM}{{\cal M}}
\newcommand{\cN}{{\cal N}}

\newcommand{\cT}{{\cal T}}

\DeclareMathOperator*{\argmin}{arg\,min}

\IEEEoverridecommandlockouts

\begin{document}
\title{Efficient Decentralized Economic Dispatch \\
for Microgrids with Wind Power Integration}

\author{\IEEEauthorblockN{Yu Zhang and Georgios B. Giannakis}
\authorblockA{Dept. of ECE and DTC, University of Minnesota, Minneapolis, USA \\
Emails: \{zhan1220, georgios\}@umn.edu}
\thanks{This work was supported by the
NSF ECCS grant 1202135, and University of Minnesota
Institute of Renewable Energy and the Environment (IREE) grant RL-0010-13.}}

\maketitle

\begin{abstract}
Decentralized energy management is of paramount importance in smart microgrids 
with renewables for various reasons including environmental friendliness, reduced communication
overhead, and resilience to failures. In this context, the present work
deals with distributed economic dispatch and demand response initiatives for
grid-connected microgrids with high-penetration of wind power. To cope
with the challenge of the wind's intrinsically stochastic availability,
a novel energy planning approach involving the actual wind energy
as well as the energy traded with the main grid, is introduced.
A stochastic optimization problem is formulated to minimize the
microgrid net cost, which includes conventional generation cost
as well as the expected transaction cost incurred by wind uncertainty.
To bypass the prohibitively high-dimensional integration involved, an efficient
sample average approximation method is utilized to obtain a
solver with guaranteed convergence. Leveraging the special infrastructure
of the microgrid, a decentralized algorithm is further developed via
the alternating direction method of multipliers. Case studies are tested
to corroborate the merits of the novel approaches.
\end{abstract}

\begin{IEEEkeywords}
Microgrids, economic dispatch, renewable energy,
sample average approximation, ADMM
\end{IEEEkeywords}

\section*{Nomenclature}

\addcontentsline{toc}{section}{Nomenclature}

\subsection{Indices, numbers, and sets}

\begin{IEEEdescription}[\IEEEusemathlabelsep\IEEEsetlabelwidth{$M$, $m$}]

\item[$T$, $t$] Number of scheduling periods, and period index.

\item[$M$, $m$] Number of conventional distributed generation (DG) units, and their index.

\item[$N$, $n$] Number of dispatchable loads, and load index.

\item[$I$, $i$] Number of wind farms, and their index.

\item[$\cM$] Set of conventional DG units.

\item[$\cN$] Set of dispatchable loads.

\end{IEEEdescription}

\subsection{Variables}

\begin{IEEEdescription}[\IEEEusemathlabelsep\IEEEsetlabelwidth{$P_{G_m}^t$}]

\item[$P_{G_m}^t$] Power output of DG unit $m$ over time slot $t$.

\item[$P_{D_n}^t$] Power consumption of load $n$ over slot $t$.

\item[$W_{i}^t$] Power output from $i$th wind farm over slot $t$.

\item[$P_R^t$] Wind power delivered to the microgrid in slot $t$.

\end{IEEEdescription}

\subsection{Constants}

\begin{IEEEdescription}[\IEEEusemathlabelsep\IEEEsetlabelwidth{$P_{D_n}^{\min,t}$, $P_{D_n}^{\max,t}$}]

\item[$P_{G_m}^{\min}$, $P_{G_m}^{\max}$] Minimum and maximum power output of conventional DG unit $m$.

\item[$R_m^{\text{up}}$, $R_m^{\text{down}}$] Ramp-up and ramp-down limits of conventional DG unit $m$.

\item[$\mathsf{SR}^t$] Spinning reserve for conventional DG.

\item[$L^t$] Fixed power demand of critical loads over slot $t$.

\item[$P_{D_n}^{\min}$, $P_{D_n}^{\max}$] Minimum and maximum power consumption of load $n$.

\item[$P_{R}^{\min}$, $P_{R}^{\max}$] Lower and upper bounds for $P_{R}^t$.

\item[$\alpha^t$, $\beta^t$] Purchase and selling prices per slot $t$.

\end{IEEEdescription}

\subsection{Functions}

\begin{IEEEdescription}[\IEEEusemathlabelsep\IEEEsetlabelwidth{$G(\cdot)$, $\hat{G}(\cdot)$}]

\item[$C_m^t(\cdot)$] Cost of conventional DG unit $m$ in slot $t$.

\item[$U_n^t(\cdot)$] Utility of load $n$ in slot $t$.

\item[$G(\cdot)$, $\hat{G}(\cdot)$] Expected and sample-averaged transaction cost across entire horizon.

\item[$\mathcal{L}_{\rho}(\cdot)$] Partially augmented Lagrangian function.

\end{IEEEdescription}

\section{Introduction}
\begin{figure}
\centering
\includegraphics[scale =0.65]{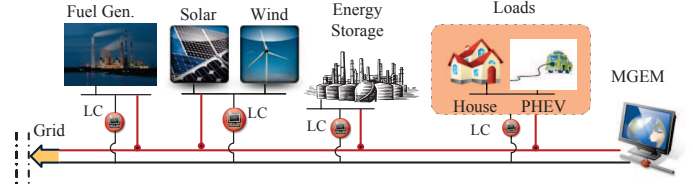}
\caption{Decentralized infrastructure of a microgrid
with communications (black) and energy flow (red) networks.
}\label{fig:MGModel}
\end{figure}

As contemporary small-scale counterparts of the bulk power grid,
smart microgrids comprise distributed energy resources (DERs) and
electricity end users, all deployed within a limited
geographical area~\cite{Hatziargyriou-PESMag}.
Depending on their origin, DERs can come either from conventional energy
sources including oil, gas and coal, or, from renewable energy sources (RES),
such as wind and solar energy. Bypassing limitations of a congested
transmission network, microgrids can generate, distribute, and regulate
power flows at the community level to efficiently meet growing consumer demands.
Besides critical non-dispatchable loads, elastic controllable
ones allow residential or commercial customers to
participate in the electricity enterprise.
A typical microgrid configuration is depicted in Fig.~\ref{fig:MGModel}.
Through the communications network, a so-termed microgrid energy manager (MGEM)
coordinates the DERs and the controllable loads, each of which has a local controller (LC).

Aligned with the goal of high-penetration RES in future smart grids,
economic dispatch (ED) with renewables has been extensively studied recently.
ED penalizing over- and under-estimation of wind energy is
investigated in~\cite{HetzerYB08}. Worst-case robust distributed
ED is proposed for grid-connected microgrids with DERs
in~\cite{YuNGGG-TSE13}. Leveraging the scenario approximation technique,
risk-constrained ED with correlated wind farms have been developed
recently in~\cite{YuNGGG-ISGT13}.
A multi-stage stochastic control approach is pursued for
risk-limiting dispatch of wind power in~\cite{RajagopalACC12}.
A chance-constrained two-stage stochastic program is formulated
in~\cite{FangGW12} for unit commitment with wind power uncertainty.
Capitalizing on the hierarchical multi-agent coordination, distributed ED
via heterogeneous wireless networks is studied in~\cite{Liang12}.

Notwithstanding their merits, the aforementioned works have limitations.
For example, it is unlikely to have the worst-case scenario
in real time operation~\cite{YuNGGG-TSE13}.
Globally optimal solutions are generally hard to obtain
for non-convex chance-constrained problems.
Convex relaxation using the scenario sampling can afford efficient optimization solvers,
but it turns out to be too conservative for scheduling the delivered renewables
in certain scenarios~\cite{YuNGGG-ISGT13}.
Moreover, slow convergence of conventional
distributed algorithms may have scalability issues facing
large-scale problems; e.g., the subgradient
ascent based dual decomposition approach~\cite{YuNGGG-TSE13}.

This paper considers day-ahead ED for microgrids with high penetration
of wind energy, operating in the grid-connected mode. By introducing
what is termed scheduled wind power, a novel energy transaction
mechanism is put forth to address the challenge of maintaining the
supply-demand balance imposed by the uncertainty of wind power.
A stochastic optimization program is formulated to
minimize the microgrid net cost, which consists of
costs for conventional generation, utility of elastic loads,
as well as the expected transaction cost (Section~\ref{sec:formulation}).
A \emph{sample average approximation (SAA)} approach with convergence guarantees is efficiently
utilized to deal with the involved multidimensional integral in the
expectation function.
With the attractive advantages of being computationally efficient and
resilient to communication outages, \emph{decentralized} scheduling over
the microgrid communications network is developed based on the
\emph{alternating direction method of multipliers (ADMM)} (Section~\ref{Section:ADMM}).
Numerical tests are reported to corroborate the merits of
the novel approaches (Section~\ref{sec:Tests}).

\noindent \emph{Notation}. Boldface lower case letters represent vectors;
$(\cdot)^{\prime}$ indicates transpose;
and $\mathbb{E}[\cdot]$ denotes the expectation operator.

\section{Robust Energy Management Formulation}\label{sec:formulation}
Consider a microgrid comprising $M$ conventional generators,
$N$ controllable (dispatchable) loads, and $I$ wind farms.
The scheduling horizon of interest is $\cT:=\{1,2,\ldots,T\}$ (e.g., one day ahead).
Let $P_{G_m}^t$ be the power produced by the $m$th conventional generator,
and $P_{D_n}^t$ the power consumed by the $n$th dispatchable load at slot $t$,
where $m \in \cM := \{1,\ldots,M\}$, $n \in \cN :=\{1,\ldots,N\}$, and $t \in \cT$.
Let $P_R^t$ denote the \emph{committed (scheduled)} wind energy 
delivered to the microgrid at slot $t$.
The ensuing subsection details the transaction mechanism between the microgrid
and the main grid. Subsection~\ref{Subsec:Form} formulates the microgrid ED problem,
which boils down to optimally dispatching the powers
$\{P_{G_m}^t\}_m$, $\{P_{D_n}^t\}_n$, and $P_R^t$ for all $t \in \cT$.

\subsection{Expected Transaction Cost}
Let $W_i^t$ denote the \emph{actual} wind power harvested from
the $i$th wind farm at time slot $t$. Suppose that the microgrid operates
in a grid-connected mode, and a transaction mechanism with the main grid
is in place, where the microgrid can buy (sell) energy
from (to) the spot market.
Specifically, the shortfall between the actual wind power produced and the one
scheduled per slot $t$ is $\left[P_{R}^t-\sum_{i=1}^IW_i^t\right]^{+}$,
while the corresponding surplus is $\left[P_{R}^t-\sum_{i=1}^IW_i^t\right]^{-}$,
where $[a]^{+}:=\max\{a,0\}$, and $[a]^{-}:= -\min\{a,0\}$.
The amount of energy shortage $\left[P_{R}^t-\sum_{i=1}^IW_i^t\right]^{+}$ is
bought with a known fixed purchase price $\alpha^t$, while the energy surplus 
$\left[P_{R}^t-\sum_{i=1}^IW_i^t\right]^{-}$ is sold back to the main grid
with a fixed selling price $\beta^t$.
Clearly, only one of these two quantities is nonzero at each slot $t$.
Wind power $W_i^t$ is a function
of the random wind speed $v_i^t$, for which different models and
wind-speed-to-wind-power mappings $W_i^t(v_i^t)$ are available~\cite{CartaRV09}.
The \emph{expected} transaction cost can be readily expressed as
\begin{align}
G(\mathbf{p}_R):=\mathbb{E}_{\mathbf{v}}\Big[\sum_{t=1}^T \Big( &\alpha^t
[P_R^t-\sum_{i=1}^{I}W_i^t(v_i^t)]^{+} \nonumber  \\
- &\beta^t[P_R^t-\sum_{i=1}^{I}W_i^t(v_i^t)]^{-}\Big)\Big]
\end{align}
where $\mathbf{v}:= [v_1^1,\ldots,v_1^T,\ldots,v_I^1,\ldots,v_I^T]^{\prime}$
and $\mathbf{p}_R := [P_R^1,\ldots,P_R^T]^{\prime}$.

\subsection{Microgrid Net Cost Minimization}
\label{Subsec:Form}

The cost of the $m$th conventional generator is a convex increasing
function $C_m^t(P_{G_m}^t)$, typically chosen either as
piecewise linear or as smooth quadratic. Moreover, the utility function of the
$n$th dispatchable load, $U_n^t(P_{D_n}^t)$, is selected to be concave increasing,
and likewise either piecewise linear or smooth quadratic.
Apart from dispatchable loads, there is also a fixed power demand from critical loads,
denoted by $L^t$. For notational brevity, let $\mathbf{p}_G$ and $\mathbf{p}_D$ denote
the vectors collecting $\{P_{G_m}^t\}_{m,t}$ and $\{P_{D_n}^t\}_{n,t}$, respectively.

ED aims at minimizing the microgrid-wide net cost; that is,
the cost of conventional generation, minus the load utility
as well as the expected transaction cost:
\begin{subequations}
\label{DERsched}
\begin{align}
\hspace{-1.7mm} \text{(P1)} \quad & \hspace{-0.3cm}\mathop{\min}_{\{\mathbf{p}_G,\mathbf{p}_D,\mathbf{p}_R\}}
\bigg\{\sum^{T}_{t=1}\left(\sum^{M}_{m=1}C_m^t(P_{G_m}^t)-\sum^{N}_{n=1}U_n^t(P_{D_n}^t)\right) \nonumber \\
&\hspace{4.5cm}+G(\{P_{R}^t\})\bigg\} \label{ObjFunc}\\
\vspace{2cm}
&\textrm{subject to:} \nonumber \\
&P_{G_m}^{\min} \le P_{G_m}^t \le P_{G_m}^{\max},~\forall~m \in \cM, ~\forall~t \in \cT\label{Plimits}\\
&P_{G_m}^t-P_{G_m}^{t-1} \le R_m^{\textrm{up}},~\forall~m \in \cM, ~\forall~t \in \cT\label{RampUp}\\
&P_{G_m}^{t-1}-P_{G_m}^t \le R_m^{\textrm{down}},~\forall~m \in \cM, ~\forall~t \in \cT\label{RampDown}\\
&\sum^{M}_{m=1}(P_{G_m}^{\max}-P_{G_m}^{t}) \ge \mathsf{SR}^{t},~\forall~t \in \cT\label{SpinRes}\\
&P_{D_n}^{\min} \le P_{D_n}^t \le P_{D_n}^{\max},~\forall~n \in \cN, ~\forall~t \in \cT\label{Dlimits}\\
&P_{R}^{\min} \le P_{R}^t \le P_{R}^{\max},~\forall~t \in \cT\label{Rlimits}\\
&\sum^{M}_{m=1}P_{G_m}^t+P_R^t = \sum^{N}_{n=1}P_{D_n}^t+L^t,~\forall~t \in \cT.\label{Balance}
\end{align}
\end{subequations}
Constraints \eqref{Plimits}-\eqref{Dlimits} stand for the minimum/maximum conventional generation,
ramping up/down limits, spinning reserves, and the minimum/maximum power
of the dispatchable loads, respectively.
They capture the typical physical limits of the power generators and loads.
Constraint~\eqref{Rlimits} places upper and lower limits on the committed wind power, which
are imposed by the capacity of the transmission lines over which the energy is transacted.
Finally, constraint~\eqref{Balance} is the supply-demand \emph{balance equation}
ensuring that the total demand is satisfied by the power generation at any time.

Note that~\eqref{Plimits}-\eqref{Balance} are linear, while $C_m^t(\cdot)$
and $-U_n^t(\cdot)$ are convex.
Consequently, the convexity of (P1) depends on that of $G(\mathbf{p}_R)$,
which is established in the following proposition.
\begin{proposition}
\label{prop:convex}
If the selling price $\beta^t$ does not exceed the purchase price $\alpha^t$ for any $t \in \mathcal{T}$,
then the expected transaction cost $G(\mathbf{p}_R)$ is convex in $\mathbf{p}_R$.
\end{proposition}
\begin{IEEEproof}
Using the identities $[a]^{+}+[a]^{-}=|a|$ and $[a]^{+}-[a]^{-}=a$,
$G(\mathbf{p}_R)$ can be equivalently re-written as
\begin{align}\label{eq:propconvex}
G(\mathbf{p}_R)=\mathbb{E}_{\mathbf{v}}\Big[\sum_{t=1}^T
\Big(&\delta^t\Big|P_R^t-\sum_{i=1}^{I}W_i^t(v_i^t)\Big| \nonumber \\
+&\gamma^t[P_R^t-\sum_{i=1}^{I}W_i^t(v_i^t)]\Big)\Big]
\end{align}
with $\delta^t := (\alpha^t-\beta^t)/2$, and $\gamma^t := (\alpha^t+\beta^t)/2$.
Since the absolute value function is convex, and the operations of nonnegative weighted summation
and integration preserve convexity~\cite[Sec.~3.2.1]{Boyd}, the claim follows readily.
\end{IEEEproof}

An immediate corollary here is that the ED problem~$(P1)$ is convex
if $\beta^t\leq \alpha^t$ for all $t$.
The next section begins with a special case when $\alpha^t \equiv \beta^t$,
before developing an approximation method together with an efficient
decentralized solver for general transaction prices satisfying the
condition of Proposition~\ref{prop:convex}.

\section{Sample Average Approximation and Distributed Algorithm}\label{Section:ADMM}
\subsection{A Special Case}

If the \emph{locational marginal pricing (LMP)} mechanism is utilized 
to price energy purchases and sales for the microgrid,
then $\alpha^t = \beta^t  = \ell^t$ for all $t \in \cT$, where $\{\ell^t\}$ are the locational marginal prices at
the bus where the transaction takes place.
In this case, we have $\delta^t = 0$ and $\gamma^t = \alpha^t$ for all $t$.
It thus follows that
\begin{align*}
G(\mathbf{p}_R)&=\mathbb{E}_{\mathbf{v}}\Big[\sum_{t=1}^T \alpha^t\Big(P_R^t-\sum_{i=1}^{I}W_i^t(v_i^t)\Big)\Big]\notag \\
& \doteq \sum_{t=1}^T \alpha^t\Big(P_R^t-\sum_{i=1}^{I}\bar{W}_i^t\Big)
\end{align*}
where $\{\bar{W}_i^t\}_{i,t}$ are sample average wind power estimates assumed to be available via
statistical inference based on historical data, or, through numerical weather prediction.

In this special case, (P1) boils down to a smooth convex minimization problem.
If $\{C_m^t(\cdot)\}_{m,t}$  and $\{U_n^t(\cdot)\}_{n,t}$
are convex quadratic or piece-wise linear, then (P1) is either a convex quadratic program (QP)
or a linear program (LP); hence, (P1) is efficiently solvable with off-the-shelf QP/LP solvers.
Next, the general case of transaction prices is investigated with the resulting optimization problem
formulated using the aforementioned sample approximation method, and solved in a decentralized fashion.

\subsection{Sample Average Approximation}
Consider now the general case under the price condition of Proposition~\ref{prop:convex},
which typically holds for microgrid power systems~\cite{Choi11}.
If the selling and buying prices are not always the same, then the absolute value terms in
$G(\mathbf{p}_R)$ do not disappear (cf.~\eqref{eq:propconvex}). Due to the nonlinearity of the absolute value operator, 
it cannot be interchanged with the expectation. 
In addition, although entries of $\mathbf{v}$ are Weibull distributed, their correlation prevents analytical 
expression of $G(\mathbf{p}_R)$.
Moreover, the multidimensional integration needed to carry out the expectation
cannot be computed with high accuracy numerically.

To bypass this challenge, the empirical estimate of $G(\mathbf{p}_R)$ will be
adopted based on $N_s$ Monte Carlo samples $\{W_i^t(s)\}_{s=1}^{N_s}$ for each $W_i^t$.
In this case, $G(\mathbf{p}_R)$ is replaced by
\begin{align}\label{eq:approxF}
\hat{G}(\mathbf{p}_R):= \frac{1}{N_s}\sum_{s=1}^{N_s}
&\sum_{t=1}^T\delta^t\Big|P_R^t-\sum_{i=1}^{I}W_i^t(s)\Big| \nonumber \\
+&\sum_{t=1}^T\gamma^t\Big(P_R^t-\sum_{i=1}^{I}\bar{W}_i^t\Big).
\end{align}
This sample average approximation (SAA) of (P1) is distribution free, and the law of large numbers (LLN)
guarantees that $\hat{G}(\mathbf{p}_R)$ is a good approximation of $G(\mathbf{p}_R)$ for $N_s$ large enough.
Based on the latter, the ED problem of interest can be approximated as
\vspace{.2cm}

\begin{subequations}\label{DERsched-app}
\hspace{-.5cm}
\fbox{
 \addtolength{\linewidth}{-2\fboxsep}%
 \addtolength{\linewidth}{-2\fboxrule}%
 \begin{minipage}{0.98\linewidth}
\begin{align}
\text{(AP1)} \mathop{\min}_{\{\mathbf{p}_G,\mathbf{p}_D,\mathbf{p}_R\}}
&\bigg\{\sum^{T}_{t=1}\Big(\sum^{M}_{m=1}C_m^t(P_{G_m}^t) \nonumber \\
&\hspace{+0.5cm}-\sum^{N}_{n=1}U_n^t(P_{D_n}^t)\Big)
+\hat{G}(\mathbf{p}_{R})\bigg\} \label{ObjFunc-AP1}\\
\textrm{s.t.}\quad &\eqref{Plimits} - \eqref{Balance}. \nonumber
\end{align}
\end{minipage}
}
\end{subequations}

\vspace{.2cm}

Clearly, convexity is preserved in the SAA formulation (AP1), and this renders it
efficiently solvable. The following conditions are sufficient to establish
the convergence of SAA applied to (P1):
A1) The optimal solution set of (P1) is nonempty;
A2) The LLN holds pointwise; that is,
$\hat{G}(\mathbf{p}_R) \rightarrow G(\mathbf{p}_R)$
with probability (w.p.) $1$ as $N_s \rightarrow \infty$.

Let $\vartheta^{*}$ and $\mathcal{S}^{*}$ denote the optimal value and the
optimal solution set of (P1), respectively.
Similarly, $\hat{\vartheta}_{N_s}$ and $\mathcal{\hat{S}}_{N_s}$ for (AP1).
Define further the deviation of the set $\mathcal{A}$ from the set $\mathcal{B}$
by $\mathbb{D(\mathcal{A},\mathcal{B})}:=
\sup_{\mathbf{x}\in \mathcal{A}}\inf_{\mathbf{y}\in \mathcal{B}}\|\mathbf{x}-\mathbf{y}\|$.
With these notational conventions, the following convergence result can be
established.

\begin{proposition}
\label{prop:converge}
If conditions A1) and A2) hold, then $\hat{\vartheta}_{N_s} \rightarrow \vartheta^{*}$,
and $\mathbb{D}(\mathcal{\hat{S}}_{N_s},\mathcal{S}^{*}) \rightarrow 0$ w.p. $1$
as $N_s \rightarrow \infty$.
\end{proposition}
\begin{IEEEproof}
It can be shown that A1)-A2) as well as the special structure of (P1)
satisfy the conditions in~\cite[Thm. 5.4]{Shapiro09},
where a convergence claim for a general problem is established.
Due to space limitations, the detailed proof is omitted.
\end{IEEEproof}

Note that (AP1) entails a separable convex objective~\eqref{ObjFunc-AP1} with
a linear equality constraint~\eqref{Balance}, as well as the compact polyhedral
feasible sets~\eqref{Plimits}-\eqref{Rlimits}, which are in the form of a Cartesian product.
This separable structure motivates solving (AP1) in a distributed fashion by resorting to
the alternating direction method of multipliers (ADMM)~\cite{BoydADMM}, which has drawn
growing interest recently, because it exhibits good performance in many large-scale distributed
optimization problems in e.g., machine learning and signal processing.

By exploiting the microgrid infrastructure, an ADMM-based distributed solver is developed
in the ensuing section.

\subsection{Decentralized ED via ADMM}

With reference to the microgrid depicted in Fig.~\ref{fig:MGModel},
it is natural to implement ED across the local controllers (LCs) of
conventional generators, dispatchable loads, and renewable facilities.
To this end, introduce a Lagrange multiplier vector
$\bm{\lambda}:=[\lambda^1,\ldots,\lambda^T]^{\prime}$
associated with the coupling equality constraints~\eqref{Balance},
along with a quadratic penalty.
The partially augmented Lagrangian of (AP1) is
\begin{align}
&\mathcal{L}_{\rho}(\mathbf{p}_G,\mathbf{p}_D,\mathbf{p}_R,\bm{\lambda})
= \sum^{T}_{t=1}\sum^{M}_{m=1}C_m^t(P_{G_m}^t)-\sum^{T}_{t=1}\sum^{N}_{n=1}U_n^t(P_{D_n}^t)\nonumber \\
&+\hat{G}(\mathbf{p}_R)
+\sum^{T}_{t=1}\lambda^t\left(\sum^{M}_{m=1}P_{G_m}^t+P_R^t-\sum^{N}_{n=1}P_{D_n}^t-L^t\right) \nonumber \\
&+\frac{\rho}{2}\sum^{T}_{t=1}\left(\sum^{M}_{m=1}P_{G_m}^t+P_R^t-\sum^{N}_{n=1}P_{D_n}^t-L^t\right)^2
\end{align}
where $\rho >0$ is a constant.

ADMM is tantamount to updating first the primal variables in the
Gauss-Seidel fashion (a.k.a. block coordinate descent), and then
updating the dual variables in a gradient ascent manner. Specifically,
with $\mathcal{P}_G := \{\mathbf{p}_G|~\eqref{Plimits}-\eqref{SpinRes}\}$,
$\mathcal{P}_D := \{\mathbf{p}_D|~\eqref{Dlimits}\}$,
and $\mathcal{P}_R := \{\mathbf{p}_R|~\eqref{Rlimits}\}$,
let $k$ denote the iteration index, and $\nu>0$ a constant stepsize.
The resulting distributed ED solver is tabulated as Algorithm~\ref{algo:Distri},
where the last step is a reasonable termination criterion using
the primal residual (see also~\cite[Sec.~3.3.1]{BoydADMM})
\begin{align}\label{eq:primres}
\xi:=\left[\sum^{T}\limits_{t=1}\left(\sum^{M}\limits_{m=1}P_{G_m}^t+P_R^t
-\sum^{N}\limits_{n=1}P_{D_n}^t-L^t\right)^2\right]^{1/2}.
\end{align}

\begin{algorithm}[t]
\caption{Distributed Economic Dispatch using ADMM}
\label{algo:Distri}
\begin{algorithmic}[1]
\State Initialize $\bm{\lambda}(0) = \mathbf{0}$
\Repeat \quad ($k = 1,2,\ldots$)
\State \textbf{Update primal variables:}
\begin{align}
\mathbf{p}_G(k+1) &= \argmin_{\mathbf{p}_G\in \mathcal{P}_G}~
\mathcal{L}_{\rho}(\mathbf{p}_G,\mathbf{p}_D(k),\mathbf{p}_R(k),\bm{\lambda}(k)) \label{eq:subPG}  \\
\mathbf{p}_D(k+1) &= \argmin_{\mathbf{p}_D\in \mathcal{P}_D}~
\mathcal{L}_{\rho}(\mathbf{p}_G(k+1),\mathbf{p}_D,\mathbf{p}_R(k),\bm{\lambda}(k)) \label{eq:subPD}  \\
\mathbf{p}_R(k+1) &= \argmin_{\mathbf{p}_R\in \mathcal{P}_R}~
\mathcal{L}_{\rho}(\mathbf{p}_G(k+1),\mathbf{p}_D(k+1),\mathbf{p}_R,\bm{\lambda}(k))  \label{eq:subPR}
\end{align}
\State \textbf{Update dual variables:} for all $t\in \mathcal{T}$
\begin{align}
\lambda^t(k+1) = \lambda^t(k) + \nu \Big(&\sum^{M}_{m=1}P_{G_m}^t(k+1)+P_R^t(k+1)  \nonumber \\
-&\sum^{N}_{n=1}P_{D_n}^t(k+1)-L^t\Big) \label{eq:subDual}
\end{align}
\Until $\xi \le \epsilon_{\textrm{res}}$
\end{algorithmic}
\end{algorithm}

\begin{figure}
\centering
\includegraphics[scale=0.9]{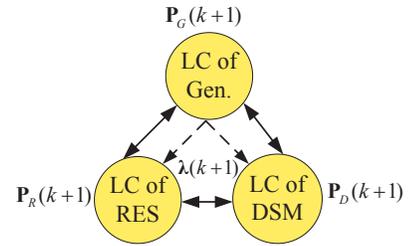}
\caption{ADMM message passing.}
\label{fig:decomp}
\end{figure}

\begin{remark}\textit{(Convergence of ADMM)}.
Sufficient conditions for convergence of the K-block
($K\ge 3$) ADMM have been established recently in~\cite{Han12} and~\cite{Mingyi12}.
One of these conditions requires that all subproblems of updating the
primal variables are strongly convex. It is worth pointing out that
although subproblem~\eqref{eq:subPR} is not strongly convex,
the algorithm always converged in the extensive numerical tests
that we performed (see Section~\ref{sec:Tests}).
Furthermore, the proximal ADMM of~\cite{Mingyi12} can be
applied here with guaranteed linear convergence.
Interested readers are referred to~\cite{Mingyi12} for
the detailed algorithm and convergence claims.
\end{remark}

ADMM iterations easily lend themselves to a distributed implementation utilizing the
microgrid communication network (cf. Fig.~\ref{fig:MGModel}).
Specifically, the LCs of conventional generation, dispatchable loads, and RES solve
subproblems~\eqref{eq:subPG}, \eqref{eq:subPD}, and~\eqref{eq:subPR} sequentially,
via efficient QP solvers.
Note that after each LC solves its own subproblem, the correspondingly updated primal variables
should be broadcast to all other LCs. The dual updating step~\eqref{eq:subDual}
can be readily implemented by any one of the three LCs.
The detailed message passing process is depicted in Fig.~\ref{fig:decomp}.

\section{Numerical Tests}
\label{sec:Tests}
In this section, case studies are presented to verify the
performance of ADMM-based distributed ED for a microgrid
consisting of $M=3$ conventional generators, $N=3$ dispatchable loads,
and $I=4$ wind farms scheduled over $T=8$ hours.
The generation costs $C_m(P_{G_m})=a_mP_{G_m}^2+b_mP_{G_m}$,
and the utilities of elastic loads $U_n(P_{D_n})=c_nP_{D_n}^2+d_nP_{D_n}$
are selected time-invariant and quadratic. The corresponding parameters of
generators, loads and transaction prices are listed in
Table~\ref{tab:gen} -- \ref{tab:price}, while spinning reserves are
set to $\mathsf{SR}^t =6.66$\,kWh for all $t \in \mathcal{T}$.
The resulting optimization problems are specified and solved via
the Matlab-based modeling language \texttt{CVX}~\cite{cvx}
along with the solver \texttt{Gurobi}~\cite{gurobi}.

To obtain the wind power samples $\{W_i^t(s)\}_{s=1}^{N_s}$ required as input to (AP1) (cf.~\eqref{ObjFunc-AP1}),
a simple but effective sampling approach leveraging autoregressive models with the wind-speed-to-wind-power mapping
is utilized; see~\cite{YuNGGG-ISGT13} and~\cite{Villanueva12} for details.
%
In the numerical tests, the sample size is $N_s=1,000$, and
the averaged wind power outputs $\{\bar{W}_i^t\}_{i,t}$ are obtained
using $20,000$ samples of the wind speed.

\begin{table}[t]
\centering
\caption{Generating capacities, ramping limits, and cost coefficients.
The units of $a_m$ and $b_m$ are \textcent/(kWh)$^{2}$
and \textcent/kWh, respectively.}\label{tab:gen}
    \begin{tabular}{  c || c | c | c | c | c }
    \hline
Unit &$P_{G_m}^{\min}$ &$P_{G_m}^{\max}$  &$R_{m,\text{up(down)}}$   &$a_m$   &$b_m$  \\ \hline
1    & 5                      & 70             & 30            &0.006         & 14   \\
2    & 5                      & 80             & 35            &0.003         & 20   \\
3    & 10                     & 85             & 50            &0.004         & 50    \\
    \hline
    \end{tabular}
\end{table}

\begin{table}[t]
\centering
\caption{load demand limits, and utility coefficients.\newline
The units of $c_n$ and $d_n$ are \textcent/(kWh)$^{2}$
and \textcent/kWh, respectively.}\label{tab:load}
    \begin{tabular}{  c || c | c | c | c }
    \hline
Unit &$P_{D_n}^{\min}$ &$P_{D_n}^{\max}$      &$c_n$       &$d_n$     \\ \hline
1    & 5                      & 30              &-0.20         & 20    \\
2    & 8                      & 50               &-0.30         & 30    \\
3    & 3                     & 45                &-0.17         & 17     \\
    \hline
    \end{tabular}
\end{table}

\begin{table}[t]
\renewcommand{\arraystretch}{1.2}
\centering
\caption{Fixed load demand and transaction prices. \newline
The unit of $\alpha^t$ and $\beta^t$ is \textcent/kWh}\label{tab:price}
    \begin{tabular}{  c || c | c | c | c | c | c | c | c }
    \hline
    \text{Time slot}    & 1    &  2   &  3    &   4  &   5  &  6    &   7   & 8 \\ \hline \hline
    $L^t$             &30    &34   &47 &60 &75 &67 &55 &43 \\ \hline
    $\alpha^t$       &1.40   &2.20   &4.70   &6.30  &8.50  &7.80   &5.60  &4.50 \\ \hline
    $\beta^t$       &1.12 &1.76  &3.76 &5.04 &6.80  &6.24 &4.48  &3.60 \\
    \hline
    \end{tabular}
\end{table}

\begin{figure}[t]
\centering
\includegraphics[width=0.45\textwidth]{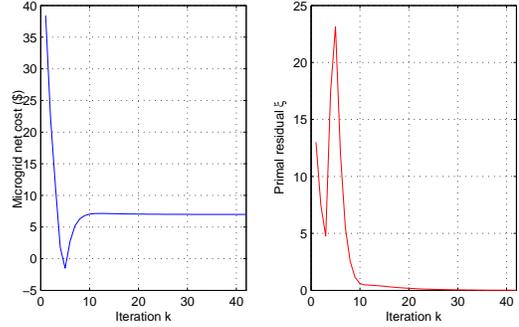}
\caption{Convergence of the net cost and evolution of the primal residual.}
\label{fig:ObjConverge}
\end{figure}

\begin{figure}[ht!]
\centering
\includegraphics[width=0.45\textwidth]{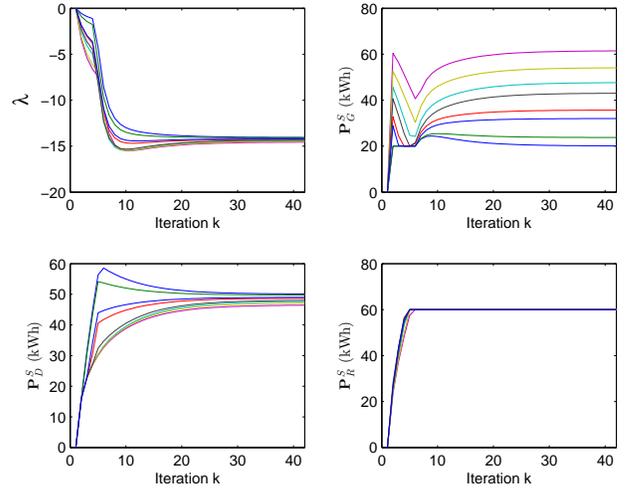}
\caption{Convergence of the primal and dual variables.}
\label{fig:VariConverge}
\end{figure}

Figure~\ref{fig:ObjConverge} demonstrates the convergence of the
net cost~\eqref{ObjFunc-AP1}, and the evolution of the primal
residual $\xi$. It is clear that the algorithm converges fast
within $50$ iterations. In all numerical tests,
the relevant parameters are $\rho=1$, $\nu=0.5$,
and $\epsilon_{\textrm{res}} = 10^{-2}$. Furthermore, as with other
distributed solvers (e.g., dual decomposition using subgradient ascent),
ADMM is not an iterative algorithm guaranteeing a monotonically
decreasing objective. Figure~\ref{fig:ObjConverge} shows that
some objective values of the iterates can be even smaller than
the optimal value due to the constraint violation.
However, for the day-ahead energy planning problem,
ADMM outperforms alternative distributed solvers thanks to its fast convergence.

Convergence of the primal and dual variables is verified in
Fig.~\ref{fig:VariConverge}, where ${P}_G^S:= \sum_{m}P_{G_m}^t$
denotes the total conventional power generation, and likewise
for $P_D^S$. Clearly, the iterates converge fast as shown by the
$8$ trajectories per subplot, each corresponding to a different
time slot.

\begin{figure}[t]
\centering
\includegraphics[width=0.45\textwidth]{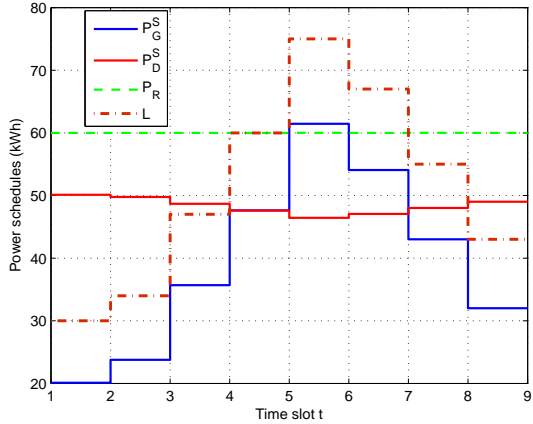}
\caption{Optimal power schedules.}
\label{fig:schedules}
\end{figure}

The optimal power schedules are depicted in Fig.~\ref{fig:schedules}.
As expected, the total conventional power generation $P_G^S$ varies
across $t$ with the same trend as the fixed load demand $L$.
Moreover, the elastic demand $P_D^S$ exhibits opposite trend with respect to $L$.
This is because when $L^t$ is low, $P_D^t$ can increase to gain in utility,
as long as the total load demand can be satisfied.
As shown in the slots from $4$ to $7$,
this behavior illustrates the peak-load shifting ability
of the proposed design. It is also interesting to see that
the optimal scheduled wind power $P_R$ is set equal to
$P_R^{\max} = 60$\,kWh across time.
This is because with the energy purchase price $\alpha^t$
being much smaller than the generation costs $\{a_m,b_m\}_{m}$
(cf.~Tables~\ref{tab:gen} and~\ref{tab:price}),
the economic scheduling decision is to reduce the
conventional generation while purchasing as much energy
as possible to keep the supply-demand balance.

\begin{figure}[t]
\centering
\includegraphics[width=0.45\textwidth]{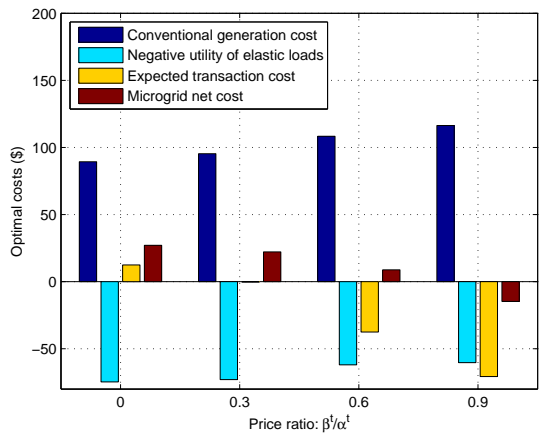}
\caption{Optimal costs with different price ratios.}
\label{fig:cost}
\end{figure}

Finally, Fig.~\ref{fig:cost} shows the effect of different transaction
prices on the optimal costs, where five times of $\alpha^t$
in Table~\ref{tab:price} is used. Clearly, the net cost decreases as
the selling-to-purchase-price ratio $\beta^t/ \alpha^t$ increases.
When this ratio increases, the microgrid can afford higher margin for revenue
by selling renewable energy back to the main grid. Thus, if more energy
is sold instead of being used within the microgrid, the cost of conventional
generation will increase to supply the loads. Therefore, as depicted in Fig.~\ref{fig:cost},
the microgrid net cost can be reduced so long as the obtained
transaction profit exceeds the extra generation cost.

\section{Conclusions and Future Work}\label{sec:Conclusions}
A distributed energy planning approach was developed in this paper tailored
for microgrids with high penetration of wind power.
By introducing the quantity of scheduled wind power, a transaction model was
proposed to maintain the supply-demand balance challenged by
the intermittent nature of RES.
A stochastic optimization problem was formulated with the objective of
minimizing the microgrid net cost. The SAA method was efficiently utilized
to overcome the high-dimensional integration involved.
Finally, the robust ED problem was solved in a distributed fashion using
an ADMM-based solver whose fast convergence was corroborated by extensive numerical tests.

A number of appealing future directions open up, including real-time dispatch and the
incorporation of uncertainty stemming from critical loads and the transaction prices.


\bibliographystyle{IEEEtran}
\bibliography{biblio}

\end{document}